\documentclass[12pt,a4paper]{article}
\usepackage[utf8]{inputenc}
\usepackage{amsmath}
\usepackage{amsfonts}
\usepackage{amssymb}
\usepackage{amsthm}
\usepackage{authblk}
\usepackage{inputenc}
\allowdisplaybreaks

\begin{document}
\title{On the Least counterexample to Robin hypothesis}
\author{Xiaolong Wu}
\affil{Ex. Institute of Mathematics, Chinese Academy of Sciences}
\affil{xwu622@comcast.net}
\date{Jun 1, 2020}
\maketitle

\begin{abstract}

    Let $G(n)=\sigma (n)/(n \log \log n )$. Robin made hypothesis that $G(n)<e^\gamma$ for all integer $n>5040$. If Robin hypothesis fails, there will be a least counterexample. This article collects the requirements the least counterexample should satisfy.
\end{abstract}
\begin{center}\textbf{ \large Introduction and notations}
\end{center}

Robin made a hypothesis [Robin 1984] that the Robin's inequality 
\begin{equation}\tag{RI}
\sigma (n)<e^\gamma n \log \log n, 
\end{equation}
 holds for all integers $n>5040$. Here $\sigma (n)=\sum_{d|n} d$ is the divisor sum function, $\gamma $ is the Euler-Mascheroni constant, log is the nature logarithm.

 For calculation convenience, we define 
 \begin{equation}\tag{1}
 \rho (n):=\frac{\sigma (n)}{n},\quad G(n):=\frac{\rho(n)}{\log \log n}.
\end{equation}
Then Robin inequality can also be written as 
\begin{equation}\tag{RI}
\rho (n)<e^\gamma \log \log n,\quad n>5040.
\end{equation}
and as 
\begin{equation}\tag{RI}
G(n)<e^\gamma, \quad n>5040.
\end{equation}

    If Robin hypothesis fails, there will be a least counterexample. This article collects the conditions the least counterexample should satisfy. Theorems A - F are known properties of the least counterexample.

    Let $n>5040$ be an integer. Write the factorization of n as 
\begin{equation}\tag{2}
n=\prod_{i=1}^{r} p_i^{a_i},
\end{equation}
where $p_i$ is the i-th prime, $p_r$ is the largest prime divisor of n.
    Theorems 1 and 2 prove that for the least counterexample n, we have
\begin{equation}\tag{3}
p_r<\log n<p_r\left(1+\frac{0.005589}{\log p_r}\right).
\end{equation}
    Theorems 3 and 4 show that for the least counterexample n, we have
\begin{equation}\tag{4}
\left\lfloor\frac{\log p_r}{\log p_i}\right\rfloor\leq a_i\leq \left\lfloor\frac{\log (a_i \log n)}{\log p_i}\right\rfloor.
\end{equation}

\noindent {\bfseries Theorem A.}
\textit{Let n be the least counterexample to Robin hypothesis. Then n is superabundant.
}
\begin{proof}
See [AF 2009] Theorem 3.
\end{proof}

\noindent {\bfseries Theorem B.}
\textit{The following properties of superabundant numbers are satisfied by the least counterexample to Robin hypothesis. Notations as in (2).\newline
1) $a_i\geq a_j$ if $i<j$.\newline
2) $\left|a_j-\left\lfloor a_j\frac{\log p_i}{\log p_j}\right\rfloor\right|\leq 1$ if $i<j$.\newline
3) $a_r=1$ except for $n=4$ and $n=36$ in chich cases we have $a_r=2$.\newline
4) $p_i^{a_i}<2^{a_1+2}$ for $2\leq i\leq r$.\newline
5) $a_i\geq \left\lfloor\frac{\log p_r}{\log p_i}\right\rfloor$ for $2\leq i\leq r$. \newline
6) Define $\epsilon (x):=\frac{1}{\log x}\left(1+\frac{3}{2\log x}\right)$ and write $\phi(n)$ for Euler totient function. Then
\[\frac{\sigma(n)}{n}>(1-\epsilon(p_r))\frac{n}{\phi(n)}.\]
}
\begin{proof}
See [Broughan 2017] Section 6.2..
\end{proof}

\noindent {\bfseries Theorem C.}
\textit{Let n, $5040<n\leq 10^{(10^{13.099} )}$, be an integer. Then n satisfies (RI). Hence the least counterexample to (RI) must be greater than $10^{(10^{13.099} )}$.
}
\begin{proof}
See [MP 2018] Theorem 5.
\end{proof}

\noindent {\bfseries Theorem D.}
\textit{Let n be the least counterexample to (RI). Define
\begin{equation}\tag{D.1}
M(k):=e^{e^{-\gamma}f(N_k)}-\log N_k,
\end{equation}
where $N_k$ is the product of first k primes and $f(n)=\prod_{p|n}\frac{p}{p-1}$. Then\newline
1) $r>969 672 728$,\newline
2) $\#\{i\leq r ; a_i\neq 1\}<\frac{r}{14}$,\newline
3) $e^{-1/\log p_r}<\frac{p_r}{\log n}<1$, \newline
4) for all $1<i\leq r$,
\begin{equation}\tag{D.2}
p_i^{a_i}<min \left(2^{a_1+2},p_i e^{M(r)}\right).
\end{equation}
}
\begin{proof}
See [Vojak 2020] Theorem 1.6.
\end{proof}

\noindent {\bfseries Theorem E.}
\textit{Let $n=\prod_{i=1}^{r} p_i^{a_i}$ be the least counterexample to (RI). Then $a_1>19$, $a_2>12$, $a_3>7$, $a_4>6$ and $a_5>5$.
}
\begin{proof}
See [Hertlein 2016] Theorems 1 and 2.
\end{proof}

\noindent {\bfseries Theorem F.}
\textit{If $n>5040$ is a sum of two squares, then n satisfies (RI).
}
\begin{proof}
See [BHMN 2008] Theorem 2.
\end{proof}

\noindent {\bfseries Theorem 1.}
\textit{Let n be the least counterexample to (RI). Then $\log n>p_r$.
}
\begin{proof}
Write $p:=p_r$. By Theorem A, we know n is superabundant, so the exponent of p in n is 1. Assume, to the contrary, $\log n\leq p$, we have
\begin{align*}
\frac{G(n)}{G(n/p)}&=\frac{\rho(n)\log\log (n/p)}{\rho(n/p)\log \log n}\\
&=\left(1+\frac{1}{p}\right)\frac{\log(\log n-\log p)}{\log \log n}\\
&=\left(1+\frac{1}{p}\right)\frac{\log\log n+\log\left(1- \frac{\log p}{\log n}\right)}{\log \log n}\\
&<\left(1+\frac{1}{p}\right)\left(1-\frac{\log p}{\log n \log\log n}\right)\\
&=1+\frac{\log n \log \log n-p\log p-\log p}{p\log n \log \log n}<1.\tag{1.1}
\end{align*}
That is, $G(n)<G(n/p)$, which means $n/p$ is also a counterexample of (RI). This contradicts to the minimality of n.
\end{proof}

\noindent {\bfseries Theorem 2. }
\textit{let $N>10^{\left(10^{13}\right)}$ be an integer. If
\begin{equation}\tag{2.1}
p_r\leq (\log N)\left(1-\frac{0.005587}{\log \log N}\right).
\end{equation}
or Conversely,
\begin{equation}\tag{2.2}
\log N > p_r\left(1+\frac{0.005589}{\log p_r}\right).
\end{equation}
Then N satisfies (RI). \newline
    Hence the least counterexample n of (RI) satisfies $\log n\leq p_r \left(1+\frac{0.005589}{\log p_r}\right)$.
}
\begin{proof}
The proof is almost identical to Theorem 9 of [Wu 2019].
\end{proof}

\noindent {\bfseries Definition 1. }
\textit{Let n be the least counterexample of (RI). Define
\begin{equation}\tag{D1.1}
x_k=(k\log n)^{1/k},\quad k=1,2,\dots,untill\quad x_k<2
\end{equation}
Then define a function
\begin{equation}\tag{D1.2}
U_n (p_i):=\left\lfloor\frac{\log(k \log n )}{\log p_i}\right\rfloor,\quad when\quad x_{k+1}<p_i\leq x_k.
\end{equation}
Since $p_r<\log n$, $U_n (p_i)$ is well defined for all primes $p_i\leq  p_r$.
}

\noindent {\bfseries Lemma 1. }
\textit{For $x_{k+1}< p_i\leq x_k$, we have $U_n(p_i)=k$.
}
\begin{proof} Since $x_{k+1}< p_i\leq x_k$, we have $\log x_{k+1}<\log p_i\leq \log x_k$. By (D1.1) 
\begin{equation}\tag{L1.1}
\frac{\log((k+1)\log n)}{k+1}<\log p_i\leq \frac{\log(k\log n)}{k}.
\end{equation}
The left inequality means
\begin{equation}\tag{L1.2}
k+1>\frac{\log((k+1)\log n)}{\log p_i}\geq \left\lfloor\frac{\log(k\log n)}{\log p_i}\right\rfloor=U_n(pi).
\end{equation}
The right inequality means
\begin{equation}\tag{L1.3}
k\leq \frac{\log(k\log n)}{\log p_i}\Longrightarrow k\leq \left\lfloor\frac{\log(k\log n)}{\log p_i}\right\rfloor=U_n (p_i ).
\end{equation}
So, we must have $k=U_n (p_i )$.
\end{proof}

\noindent {\bfseries Theorem 3. }
\textit{Let $n=\prod_{i=1}^r p_i^{a_i} >10^{(10^13 )}$ be an integer. Assume $p_s>U_n(p_s)$ for some index s. Then $G(n)<G(n/p_s)$. \newline
This means that if n is the least counterexample of (RI), then $a_i\leq U_n(p_i)$ for all i, $1\leq i\leq r$.
}
\begin{proof} By definition of $U_n$, we have $x_{k+1}<p_s\leq x_k$, for some k. $a_s>U_n(p_s)$ means that 
\begin{equation}\tag{3.1}
\frac{\log((k+1)\log n)}{k+1}<\log p_s\leq \frac{\log(k\log n)}{k} 
\end{equation}
and $a_s>U_n(p_s)=k$. Hence $a_s\geq k+1$, and
\begin{equation}\tag{3.2}
\log p_s>\frac{\log((k+1)\log n)}{k+1}\geq \frac{\log(a_s\log n)}{a_s}.
\end{equation}
We have $p_s^{a_s}>a_s\log n$, and hence 
\begin{equation}\tag{3.3}
p_s>(a_s\log n)^{1/a_s}.
\end{equation}
Write $n_1=n/p_s$. It is easy to verify that
\begin{align*}
\frac{G(n)}{G(n_1)}&=\frac{\rho(n)\log \log n_1}{\rho(n_1)\log\log n}\\
&<\left(1+\frac{1}{p_s^{a_s}}\right)\left(1-\frac{\log p_s}{\log n \log\log n}\right)\\
&\leq \left(1+\frac{1}{a_s\log n}\right)\left(1-\frac{\log(a_s\log n )}{a_s\log n\log\log n}\right)\\
&=1+\frac{1}{a_s\log n} -\frac{\log\log n+\log a_s}{a_s\log n\log\log n}-\frac{\log(a_s\log n)}{(a_s\log n)^2\log\log n}\\
&=1-\frac{\log a_s}{a_s\log n\log\log n}-\frac{\log(a_s\log n)}{(a_s\log n)^2\log\log n}<1.\tag{3.4}
\end{align*}
That is, $G(n)<G(n/p_s )$.
\end{proof}

\noindent {\bfseries Definition 2. }
\textit{Define
\begin{equation}\tag{D2.1}
L(p_i)=L_{p_r}(p_i):=\left\lfloor \frac{\log p_r}{\log p_i}\right\rfloor \quad for i\leq r.
\end{equation}
}

\noindent {\bfseries Theorem 4.}
\textit{Let $n>10^{(10^{13} )}$ be an integer. If $a_s<L(p_s)$ for some index $s<r$, then $G(n)<G(np_s/p_r)$.\newline
Hence the least counterexample of (RI) must have $a_i\geq L(p_i)$ for all i, $1\leq i\leq r$.
}
\begin{proof}As n being superabundant, we know $a_r=1=L(p_r)$. Define $n_1=\frac{p_s}{p_r}n$. Then $n_1<n$. $a_s<L(p_s )=\left\lfloor\frac{\log p_r}{\log p_s}\right\rfloor$ means $a_s+1\leq \left\lfloor\frac{\log p_r}{\log p_s}\right\rfloor\leq \frac{\log p_r}{\log p_s}$. Hence $p_s^{a_s+1}\leq p_r$ and 
\begin{equation}\tag{4.1}
\log p_s \leq \frac{1}{a_s+1}\log p_r.
\end{equation}
It is easy to deduce
\begin{align*}
\frac{G(n)}{G(n_1)}&=\frac{\rho(n)\log \log n_1}{\rho(n_1)\log\log n}\\
&=\left(1+\frac{1}{p_r}\right)\left(\frac{p_s(p_s^{-a_s}+\dots+1)}{p_s^{a_s+1}+\dots+1}\right)\left(\frac{\log(\log n-\log p_r+\log p_s)}{\log\log n}\right)\\
&<\left(1+\frac{1}{p_r}\right)\left(1-\frac{1}{p_s^{a_s+1}+\dots+1}\right)\left(1-\frac{\log p_r-\frac{1}{a_s+1}\log p_r}{\log n \log\log n}\right)\\
&=\left(1+\frac{1}{p_r}\right)\left(1-\frac{1}{p_s^{a_s+1}+\dots+1}\right)\left(1-\left(\frac{a_s}{a_s+1}\right)\frac{\log p_r}{\log n \log\log n}\right)\tag{4.2}
\end{align*}
By Theorem 2,
\begin{equation}\notag
\log n\leq p_r\left(1+\frac{0.005589}{\log p_r}\right).
\end{equation}
Noting $n>10^{(10^{13})}$, we have
\begin{equation}\tag{4.3}
\log n\leq cp_r,\quad\quad c:=1+\frac{0.005589}{\log (2.3\times 10^{13})}=1.000235.
\end{equation}
Since $\log(cp_r)<c\log p_r$, (4.2) can be simplified to
\begin{align*}
\frac{G(n)}{G(n_1)}&<\left(1+\frac{1}{p_r}\right)\left(1-\frac{1}{p_s^{a_s+1}+\dots+1}\right)\left(1-\left(\frac{a_s}{a_s+1}\right)\frac{\log p_r}{(cp_r)\log (cp_r)}\right)\\
&<\left(1+\frac{1}{p_r}\right)\left(1-\frac{1}{p_s^{a_s+1}+\dots+1}\right)\left(1-\left(\frac{a_s}{a_s+1}\right)\frac{1}{c^2p_r}\right).\tag{4.4}
\end{align*}
Now we split the proof into two cases.\newline
Case 1) $a_s=1$. We have in this case
\begin{equation}\tag{4.5}
1-\left(\frac{a_s}{a_s+1}\right)\frac{1}{c^2 p_r}<1-\frac{1}{2c^2 p_r}<1-\frac{0.49}{p_r}.
\end{equation}
\begin{equation}\tag{4.6}
1-\frac{1}{p_s^2+p_s+1}\leq 1-\frac{4}{7p_s^{a_s+1}}<1-\frac{0.57}{p_r}.
\end{equation}
Substitute (4.5) and (4.6) in to (4.4), we get 
\begin{equation}\tag{4.7}
\frac{G(n)}{G(n_1)}<\left(1-\frac{0.49}{p_r}\right)\left(1-\frac{0.57}{p_r}\right)\left(1+\frac{1}{p_r}\right)<1.
\end{equation}
Hence $G(n)<G(n_1)$.\newline
Case 2) $a_s>1$. We have
\begin{equation}\tag{4.8}
1-\left(\frac{a_s}{a_s+1}\right)\frac{1}{c^2 p_r}<1-\frac{2}{3c^2 p_r}<1-\frac{0.66}{p_r}.
\end{equation}
\begin{equation}\tag{4.9}
1-\frac{1}{p_s^{a_s+1}+\dots+1}\leq 1-\frac{1}{2p_s^{a_s+1}}<1-\frac{0.50}{p_r}.
\end{equation}
Substitute (4.8) and (4.9) into (4.4), we get
\begin{equation}\tag{4.10}
\frac{G(n)}{G(n_1)}<\left(1-\frac{0.66}{p_r}\right)\left(1-\frac{0.50}{p_r}\right)\left(1+\frac{1}{p_r}\right)<1.
\end{equation}
Hence $G(n)<G(n_1)$.
\end{proof}

\noindent {\bfseries Theorem 5. }
\textit{
Let $n=\prod_{i=1}^r p_i^{a_i}>10^{(10^{13})}$ be the least counterexample of (RI). Let s be the largest index such that $a_s\neq 1$. Then $p_s<1.414342\sqrt{p_r}$.
}
\begin{proof} Since we are searching for the largest index, we may assume $a_s=2$ and set $k=2$ in the definition of $U_n(p_s)$. By Theorem 3 we have
\begin{equation}\tag{5.1}
2=a_s\leq U_n(p_s)=\left\lfloor \frac{\log(2\log n)}{\log p_s}\right\rfloor \leq \frac{\log(2\log n)}{\log p_s}.
\end{equation}
By Theorem 2,
\begin{equation}\tag{5.2}
p_s^2\leq 2\log n\leq 2p_r \left(1+\frac{0.005589}{\log p_r}\right).
\end{equation}
\end{proof}

\noindent {\bfseries Theorem 6. }
\textit{
Let $n=\prod_{i=1}^r p_i^{a_i}>10^{(10^{13})}$ be the least counterexample of (RI). Let s be the largest index such that $a_s\neq 1$. Then $p_s>0.999999\sqrt{p_r}$.
}
\begin{proof} Let integer t be the index such that $p_t$ is the prime just below $\sqrt{p_r}$. Then by theorem 9, we have $a_t\geq L(p_t)=\left\lfloor\frac{\log p_r}{\log p_t}\right\rfloor=2$. By Corollary 5.5 of [Dusart 2018], for all $x\geq 468 991 632$, there exists a prime p such that
\begin{equation}\tag{6.1}
x<p\leq x\left(1+\frac{1/5000}{(\log x)^2}\right).
\end{equation}
Hence $p_s\geq p_t>\left(1-\frac{1/5000}{\log \sqrt{p_r})^2}\right)\sqrt{p_r}>0.999999\sqrt{p_r}$. 
\end{proof}

\begin{center}
{\bfseries \large References}
\end{center}

\noindent {[AF 2009{]} A. Akbary and Z. Friggstad. \textit{Superabundant numbers and the Riemann hypothesis}. Amer. Math. Monthly, 116(3):273275, 2009.\\
\noindent {[}BHMN 2008{]} W. Banks, D. Hart, P Moree and W. Nevans. \textit{The Nicolas and Robin inequalities with sums of two squares}. Monatsh Math (2009) 157:303–322.\\
\noindent {[}Briggs 2006{]} K. Briggs. \textit{Abundant numbers and the Riemann hypothesis}. Experiment. Math., 15(2):251–256, 2006.\\
{[}Broughan 2017{]} K. Broughan, \textit{Equivalents of the Riemann Hypothesis} Vol 1. Cambridge Univ. Press. (2017)\\
{[}Dusart 2018{]} P. Dusart. \textit{Explicit estimates of some functions over primes}. Ramanujan J., 45(1):227–251, 2018.\\
{[}Hertlein 2016{]} A. Hertlein. \textit{ROBIN’S INEQUALITY FOR NEW FAMILIES OF INTEGERS} 2016-12-15\\
https://arxiv.org/abs/1612.05186\\
{[}Morrill;Platt 2018{]} T. Morrill, D. Platt. \textit{Robin’s inequality for 25-free integers and obstacles to analytic improvement} \\
https://arxiv.org/abs/1809.10813\\
{[}Vojak 2020{]} R. Vojak. \textit{On numbers satisfying Robin’s inequality, properties of the next counterexample and improved specific bounds}. 2020-05-19\\
https://arxiv.org/abs/2005.09307\\
{[}Wu 2019{]} X. Wu. \textit{Properties of counterexample to Robin hypothesis}. 2019-01-16\\
https://arxiv.org/abs/1901.09832
\end{document}